# Spectral Characterization of Hyponormal Weighted Shifts [*]

WANG Gong-bao[1], MA Ji-pu[2]

*(1.Dept. of Basic Courses, Naval University of Engineering, Hubei Wuhan 430033, China; 2. Dept. of Math., Nanjing University, Jiangsu Nanjing 210093, China)*

**Abstract:** The complete characterizations of the spectra and their various parts of hyponormal unilateral and bilateral weighted shifts are presented respectively in this paper. The results obtained here generalize the corresponding work of the references.

**MR(2000) Subject Classification:** 47A10, 47B20

## §1.  Introduction

Let $H$ be a separable complex Hilbert space with an orthonormal basis $\{e_n\}_{n=1}^{\infty}$, and $B(H)$ the Banach algebra of all bounded linear operators on $H$.

Suppose $S \in B(H)$. If there is a bounded sequence of complex numbers $\{\alpha_n\}_{n=1}^{\infty}$ such that $Se_n = \alpha_n e_{n+1}$ for all $n \geq 1$, then $S$ is called a unilateral weighted shift with weight sequence $\{\alpha_n\}_{n=1}^{\infty}$.

It is known that weighted shifts play an important role in the theory of Hilbert spaces. In this paper, we will give respectively the complete spectral characterizations of hyponormal unilateral and bilateral weighted shifts. And the results obtained here generalize the corresponding work of the references [1,2].

Throughout this paper, $\sigma(.), \sigma_p(.), \sigma_c(.), \sigma_r(.), \Pi(.)$ and $\Gamma(.)$ denote the spectrum, point spectrum, continuous spectrum, residual spectrum, approximate point spectrum and compression spectrum of a given operator respectively([1],[4]).

## §2.  Spectral Characterizations of Hyponormal Unilateral Weighted Shifts

Let $S$ be a unilateral weighted shift with weight sequence $\{\alpha_n\}_{n=1}^{\infty}$. We know that $S$ is hyponormal if and only if $|\alpha_n| \leq |\alpha_{n+1}|$ for all $n \geq 1$. For a hyponormal unilateral weighted shift $S$, we may assume $\alpha_n \neq 0$ for all $n \geq 1$. In fact, if $\alpha_1 = \alpha_2 = \cdots = \alpha_m = 0$ for some $m \geq 1$, then $M=\mathrm{span}\{e_1, e_2, \cdots, e_m\}$ reduces $S$ and $S$ is normal on $M$ with $S|_M = 0$. In this case, it is sufficient to study $S|_{M^\perp}$. For this reason, without real loss of generality, we assume all weight sequences to have no zeros throughout this section.

Moreover, for the above operator $S$, it is easy to see that $S^* e_1 = 0$ and $S^* e_n = \overline{\alpha_{n-1}} e_{n-1}$ for all $n \geq 2$, where $S^*$ denotes the self-adjoint operator of $S$.

Now, the main results of this section are as follows.

**Theorem 1** Let $S$ be a hyponormal unilateral weighted shift with weight sequence $\{\alpha_n\}_{n=1}^\infty$, and let $M = \lim_{n \to +\infty} |\alpha_n|$. Then

(1) $\|S\| = M$;

(2) $\sigma_p(S) = \phi$, $\sigma_p(S^*) = \{\lambda \in C \ \| \lambda \| < M\}$;

(3) $\sigma_r(S) = \{\lambda \in C \ \| \lambda \| < M\}$, $\sigma_c(S) = \{\lambda \in C \ \| \lambda \| = M\}$, $\sigma(S) = \sigma(S^*) = \{\lambda \in C \ \| \lambda \| \leq M\}$;

(4) $\sigma_r(S^*) = \phi$, $\sigma_c(S^*) = \{\lambda \in C \ \| \lambda \| = M\}$, $\Gamma(S^*) = \phi$, $\Pi(S^*) = \{\lambda \in C \ \| \lambda \| \leq M\}$;

(5) $\Gamma(S) = \{\lambda \in C \ \| \lambda \| < M\}$, $\Pi(S) = \{\lambda \in C \ \| \lambda \| = M\}$.

**Proof.** (1) It is easy to prove that $\|S\| = M$. So $\sigma(S) \subset \{\lambda \in C \ \| \lambda \| \leq M\}$.

(2) For every $\lambda \in C$, we need to prove that $N(\lambda I - S) = \{0\}$. Now, suppose there exists some $x \in H$ such that $(\lambda I - S)x = 0$. Then $Sx = \sum_{n=1}^\infty (Sx, e_n)e_n = \sum_{n=2}^\infty \alpha_{n-1}(x, e_{n-1})e_n = \sum_{n=1}^\infty \lambda(x, e_n)e_n$.

Hence $\lambda(x, e_1) = 0$ and $\alpha_{n-1}(x, e_{n-1}) = \lambda(x, e_n)$ for all $n \geq 2$. Since $\alpha_n \neq 0$ for all $n \geq 1$, so $x = 0$. This is to say that $N(\lambda I - S) = \{0\}$. Therefore, $\sigma_p(S) = \phi$.

We want to prove that $\sigma_p(S^*) = \{\lambda \in C \ \| \lambda \| < M\}$. Suppose that $0 < |\lambda_0| < M$. Set $\beta_1 = 1$ and $\beta_{n+1} = \dfrac{\beta_n \lambda_0}{\alpha_n}$ for all $n \geq 1$. Since $\lim_{n \to +\infty} \left|\dfrac{\beta_{n+1}}{\beta_n}\right| = \lim_{n \to +\infty} \left|\dfrac{\lambda_0}{\alpha_n}\right| = \dfrac{|\lambda_0|}{M} < 1$, hence

$x_0 = \sum_{n=1}^\infty \beta_n e_n \in H (x_0 \neq 0)$. It follows that $(\lambda_0 I - S^*)x_0 = 0$ and $\lambda_0 \in \sigma_p(S^*)$. It is easy to see that $0 \in \sigma_p(S^*)$. So $\{\lambda \in C \ \| \lambda \| < M\} \subset \sigma_p(S^*)$.

On the other hand, for $|\lambda| = M$, we can show that $N(\lambda I - S^*) = \{0\}$. Suppose $x \in H$ and $S^*x = \lambda x$. Then we have $\sum_{n=1}^{\infty}(S^*x, e_n)e_n = \sum_{n=1}^{N}(x, \alpha_n e_{n+1})e_n = \sum_{n=1}^{\infty}\lambda(x, e_n)e_n$ and

$$\sum_{n=1}^{\infty}|(x, \alpha_n e_{n+1})|^2 = |\lambda|^2 \|x\|^2 = M^2\|x\|^2.$$

For $\sum_{n=1}^{\infty}|(x, \alpha_n e_{n+1})|^2 \leq M^2 \sum_{n=1}^{\infty}|(x, e_{n+1})|^2 = M^2(\|x\|^2 - |(x, e_1)|^2)$, so $(x, e_1) = 0$. From the fact that $(x, \alpha_n e_{n+1}) = \lambda(x, e_n)$ and $\alpha_n \neq 0$ for all $n \geq 1$, we have $(x, e_n) = 0$ for all $n \geq 1$. Hence $x=0$ and $\sigma_p(S^*) \subset \{\lambda \in C \mid |\lambda| < M\}$. Therefore, (2) holds.

(3) For every $\lambda \in C$, we have $\overline{R(\lambda I - S)} = (N(\overline{\lambda} I - S^*))^{\perp}$. Hence $\overline{R(\lambda I - S)} \neq H$ if and only if $\overline{\lambda} \in \sigma_p(S^*)$. From (2) it follows that $\sigma_r(S) = \overline{\sigma_p(S^*)} = \{\lambda \in C \mid |\lambda| < M\}$. Since $\sigma_r(S) \subset \sigma(S) \subset \{\lambda \in C \mid |\lambda| \leq M\}$ and $\sigma(S)$ is a closed set, so $\sigma(S) = \sigma(S^*) = \{\lambda \in C \mid |\lambda| \leq M\}$. For $\sigma_p(S) = \phi$, hence $\sigma_c(S) = \{\lambda \in C \mid |\lambda| = M\}$.

(4) For every $\lambda \in C$, we have $\overline{R(\lambda I - S^*)} = (N(\overline{\lambda} I - S))^{\perp} = H$. It is evident that $\sigma_r(S^*) = \Pi(S^*) = \phi$. From (2) and (3), it follows that $\sigma_c(S^*) = \{\lambda \in C \mid |\lambda| = M\}$ and $\Pi(S^*) = \{\lambda \in C \mid |\lambda| \leq M\}$.

(5) Since $\sigma_r(S) = \Pi(S) - \sigma_p(S)$ and $\sigma_p(S) = \phi$, so $\Pi(S) = \sigma_r(S) = \{\lambda \in C \mid |\lambda| < M\}$.

In order to characterize $\Pi(S)$, the following concept and notation are needed([3]).

For any $T \in B(H)$, $m(T) = \inf\{\|Tx\| \mid \|x\| = 1\}$ is defined as the lower bound of $T$. Therefore, for every $x \in H$, we have $\|Tx\| \geq m(T)\|x\|$. From [3] we know that $\lim_{n \to +\infty}[m(T^n)]^{\frac{1}{n}}$ exists. Let

$$r_1(T) = \lim_{n \to +\infty}[m(T^n)]^{\frac{1}{n}}.$$

Particularly, for the above operator $S$, we have that $m(S^n) = \inf_{k \geq 1}|\alpha_k \alpha_{k+1} \cdots \alpha_{k+n-1}| = |\alpha_1 \alpha_2 \cdots \alpha_n|$ for all $n \geq 1$. Hence $r_1(S) = \lim_{n \to +\infty}|\alpha_1 \alpha_2 \cdots \alpha_n|^{\frac{1}{n}} = \lim_{n \to +\infty}\frac{|\alpha_1| + |\alpha_2| + \cdots + |\alpha_n|}{n} = \lim_{n \to +\infty}|\alpha_n| = M$. From Proposition 13 of

[3], it follows that if $\lambda \in \Pi(S)$, then $|\lambda| \geq r_1(S) = M$. Since $\partial\sigma(S) \subset \Pi(S)$, where $\partial\sigma(S)$ is the boundary of $\sigma(S)$, therefore $\Pi(S) = \{\lambda \in C \mid |\lambda| = M\}$.

We have completed the proof of Theorem 1.

The following corollary holds from the above theorem.

**Corollary 1** Every closed disk of the complex plane with center at the origin is the spectrum of some hyponormal unilateral weighted shift. And every closed disk is the spectrum of some hyponormal operator of the form $S + \alpha I$, where $S$ is a hyponormal unilateral weighted shift and $\alpha \in C$.

### § 3. Spectral Characterizations of Hyponormal Bilateral Weighted Shifts

In this section, we will present the complete spectral characterizations of hyponormal bilateral weighted shifts.

Let $H$ be a separable complex Hilbert space with an orthonormal basis $\{e_n\}_{n=-\infty}^{+\infty}$, and $T \in B(H)$. We call $T$ a bilateral weighted shift if there exists a bounded sequence of complex numbers $\{\beta_n\}_{n=-\infty}^{+\infty}$ such that $Te_n = \beta_n e_{n+1}$ for all $n = 0, \pm 1, \pm 2, \cdots$. In this case, $\{\beta_n\}_{n=-\infty}^{+\infty}$ is called the weight sequence of $T$. It is easy to see that $T$ is hyponormal if and only if $|\beta_n| \leq |\beta_{n+1}|$ for all $n = 0, \pm 1, \pm 2, \cdots$. And it is evident that $T^* e_n = \overline{\beta_{n-1}} e_{n-1}$ for all $n = 0, \pm 1, \pm 2, \cdots$.

Similarly, for hyponormal bilateral weighted shifts, we assume all weight sequences to have no zeros in the following discussion.

**Theorem 2** Let $T$ be a hyponormal bilateral weighted shift with weight sequence $\{\beta_n\}_{n=-\infty}^{+\infty}$, and let $M = \lim_{n \to +\infty} |\beta_n|$, $m = \lim_{n \to -\infty} |\beta_n|$. Then

(1) $\|T\| = M$;

(2) $\sigma_p(T) = \phi$;

(3) $\sigma_p(T^*) = \{\lambda \in C \mid m < |\lambda| < M\}$;

(4) $\sigma_r(T) = \{\lambda \in C \mid m < |\lambda| < M\}$, $\sigma_c(T) = \{\lambda \in C \mid |\lambda| = m, \text{or} M\}$, $\sigma(T) = \sigma(T^*) = \{\lambda \in C \mid m \leq |\lambda| \leq M\}$;

(5) $\sigma_r(T^*) = \phi$, $\sigma_c(T^*) = \{\lambda \in C \mid |\lambda| = m, \text{or} M\}$, $\Gamma(T^*) = \phi$, $\Pi(T^*) = \{\lambda \in C \mid m \leq |\lambda| \leq M\}$;

(6) $\Gamma(T) = \{\lambda \in C \mid m < |\lambda| < M\}$, $\Pi(T) = \{\lambda \in C \mid |\lambda| = m, \text{or} M\}$.

**Proof.** (1) It is easy to show that $\|T\| = M$. Hence $\sigma(T) \subset \{\lambda \in C \mid |\lambda| \leq M\}$.

(2) For every $\lambda \in C$, we want to prove that $N(\lambda I - T)=\{0\}$. Now, suppose there exists some $x \in H$ such that $(\lambda I - T)x = 0$. Then

$$Tx = \sum_{n=-\infty}^{+\infty}(Tx, e_n)e_n = \sum_{n=-\infty}^{+\infty}\beta_{n-1}(x, e_{n-1})e_n = \sum_{n=-\infty}^{+\infty}\lambda(x, e_n)e_n.$$ Hence $\lambda(x, e_n) = \beta_{n-1}(x, e_{n-1})$ for all $n = 0, \pm 1, \pm 2, \cdots$. When $\lambda = 0$, since $\beta_n \neq 0$ for all $n$, so $x=0$.

Suppose $\lambda \neq 0$ and $|\lambda| \leq M$. From $\lambda(x, e_n) = \beta_{n-1}(x, e_{n-1})$, we have the following formulas:

$$(x, e_n) = \frac{\beta_0 \beta_1 \cdots \beta_{n-1}}{\lambda^n}(x, e_0) \quad \text{for all } n \geq 1, \tag{$1^0$}$$

and $(x, e_{-n}) = \dfrac{\lambda^n}{\beta_{-1}\beta_{-2}\cdots\beta_{-n}}(x, e_0)$ for all $n \geq 1$. $\tag{$2^0$}$

In ($1^0$), let $a_n = \dfrac{\beta_0\beta_1\cdots\beta_{n-1}}{\lambda^n}$. Then $\lim\limits_{n\to+\infty}\dfrac{|a_{n+1}|}{|a_n|} = \lim\limits_{n\to+\infty}\dfrac{|\beta_n|}{|\lambda|} = \dfrac{M}{|\lambda|}$. Hence, when $|\lambda| < M$, $\{|a_n|\}$ is an increasing sequence. So $\lim\limits_{n\to+\infty}|a_n| > 0$. Since $\lim\limits_{n\to+\infty}|(x, e_n)| = \lim\limits_{n\to+\infty}|a_n|\cdot|(x, e_0)| = 0$, so $(x, e_0) = 0$. Now, from ($1^0$) and ($2^0$), we have $(x, e_n) = 0$ for all $n = 0, \pm 1, \pm 2, \cdots$. This means $x = 0$.

Suppose $|\lambda| = M$. From ($2^0$) we have

$$|(x, e_0)| = \frac{|\beta_{-1}\beta_{-2}\cdots\beta_{-n}|}{|\lambda|^n}|(x, e_{-n})| \leq \frac{M^n}{|\lambda|^n}|(x, e_{-n})| = |(x, e_{-n})|.$$ Since $\lim\limits_{n\to+\infty}|(x, e_{-n})| = 0$, so $(x, e_0) = 0$. It is easy to see that $x = 0$ from ($1^0$) and ($2^0$). We have showed that for every $\lambda \in C$, $N(\lambda I - T)=\{0\}$. Hence $\sigma_p(T) = \phi$.

(3) Suppose $\lambda \in C$ and $x \in H$ such that $T^*x = \lambda x$. Similarly, we have $\lambda(x, e_n) = \overline{\beta_n}(x, e_{n+1})$ for all $n = 0, \pm 1, \pm 2, \cdots$. And it is not difficult to show that if $|\lambda| \leq m$ or $|\lambda| = M$, then $x=0$. In this case, $\lambda \notin \sigma_p(T^*)$.

Suppose $m < |\lambda| < M$. Choose $b_0 = c_0 = 1$ and let $b_{n+1} = \dfrac{b_n\lambda}{\overline{\beta_n}}$ for all $n=0,1,2,\cdots$, and

$c_n = \dfrac{\overline{\beta_n} c_{n+1}}{\lambda}$ for all $n=-1,-2,\cdots$. It is easy to show that $x_0 = \sum_{n=-\infty}^{-1} c_n e_n + \sum_{n=0}^{+\infty} b_n e_n \in H (x_0 \neq 0)$ such that $T^* x_0 = \lambda x_0$. In this case, $\lambda \in \sigma_p(T^*)$. Thus, we conclude that $\sigma_p(T^*) = \{\lambda \in C \mid m < |\lambda| < M\}$.

(4) For every $\lambda \in C$, we have $\overline{R(\lambda I - T)} = (N(\overline{\lambda} I - T^*))^\perp$. Hence, $\sigma_r(T) = \sigma_p(T^*) = \{\lambda \in C \mid m < |\lambda| < M\}$ follows directly from (2) and (3).

Suppose $m=0$. Since the spectrum of an operator is a compact set, it is easy to prove that $\sigma(T) = \sigma(T^*) = \{\lambda \in C \mid |\lambda| \leq M\}$ and $\sigma_c(T) = \{\lambda \in C \mid \|\lambda\| = 0, or M\}$.

Let $m>0$ and $|\lambda| < m$. We turn to prove $\lambda \in \rho(T)$, where $\rho(T)$ denotes the regular set of $T$. In fact, for every $x \in H$, we have

$\|(\lambda I - T)x\| \geq \|Tx\| - |\lambda|\|x\| = (\sum_{n=-\infty}^{+\infty} |\beta_n|^2 |(x,e_n)|^2)^{\frac{1}{2}} - |\lambda|\|x\| \geq (m - |\lambda|)\|x\|$. So $\lambda I - T$ is bounded below by $m - |\lambda|$. Moreover, $\overline{R(\lambda I - T)} = (N(\overline{\lambda} I - T^*))^\perp = H$ follows from (3). Hence $\lambda \in \rho(T)$ ([4]). Because of the compactness of spectrum, we know that $\sigma(T) = \sigma(T^*) = \{\lambda \in C \mid m \leq |\lambda| \leq M\}$ and $\sigma_c(T) = \{\lambda \in C \mid |\lambda| = m, or M\}$.

(5) For every $\lambda \in C$, since $\overline{R(\lambda I - T^*)} = (N(\overline{\lambda} I - T))^\perp = H$, so $\sigma_r(T^*) = \Gamma(T^*) = \phi$. From (3) and (4), we have $\sigma_c(T^*) = \{\lambda \in C \mid |\lambda| = m, or M\}$ and $\Pi(T^*) = \{\lambda \in C \mid m \leq |\lambda| \leq M\}$.

(6) Since $\sigma_r(T) = \Gamma(T) - \sigma_p(T)$, hence $\Gamma(T) = \{\lambda \in C \mid m < |\lambda| < M\}$ follows from (2) and (4). Because the weight sequence of $T$ has no zeros, so $T$ is injective. From [2] and Proposition 15 of [3], it is true that $\Pi(T) = \{\lambda \in C \mid |\lambda| = m, or M\}$. The proof of Theorem 2 is complete.

Note that, if $\beta_n = 1$ for all $n = 0, \pm 1, \pm 2, \cdots$, then Theorem 2 includes the answer to Problem 68 of [1].

By Theorem 2, the following corollary holds.

**Corollary 2** Every closed annulus of the complex plane with center at the origin is the spectrum of some hyponormal bilateral weighted shift. And every closed annulus is the spectrum of some hyponormal operator of the form $T + \beta I$, where $T$ is a hyponormal bilateral weighted shift and $\beta \in C$.

---

\* **Research supported by the National Natural Science Foundation of China(grant no.19971039) and the Science Foundation of Naval University of Engineering(2002706).**